\documentclass[11pt]{amsart}
\usepackage[small,width=\linewidth]{caption}

\def\cput(#1,#2)#3{\put(#1,#2){\hbox to 0pt{\hss{#3}\hss}}}
\def\lput(#1,#2)#3{\put(#1,#2){\hbox to 0pt{\hss{#3}}}}

\title[\resizebox{5.3in}{!}{Construction of breather solutions for nonlinear Klein-Gordon equations on periodic metric graphs}]{Construction of breather solutions for nonlinear Klein-Gordon equations on periodic metric graphs}

\author{Daniela Maier}

\address{Daniela Maier (daniela.maier@mathematik.uni-stuttgart.de), Institut f\"{u}r Analysis, Dynamik und Modellierung,
Universit\"{a}t Stuttgart \\
Pfaffenwaldring 57, D-70569 Stuttgart, Germany}

\usepackage{amsmath,amsfonts,amssymb,amsthm,overpic,relsize}
\usepackage{mathtools}
\usepackage[square,numbers]{natbib}
\usepackage[a4paper]{geometry}
\usepackage{float}

\newtheorem{thm}{Theorem}[section]

\newtheorem{rem}[thm]{Remark}

\newtheorem{hyp}[thm]{Hypothesis}

\newcommand{\R}{{\mathbb R}}
\newcommand{\Q}{{\mathbb Q}}
\newcommand{\N}{{\mathbb N}}
\newcommand{\Z}{{\mathbb Z}}
\newcommand{\C}{{\mathbb C}}
\newcommand{\fish}{$\rtimes$}
\DeclarePairedDelimiter{\abs}{\lvert}{\rvert}
\DeclarePairedDelimiter{\norm}{\lVert}{\rVert}

\begin{document}
\long\def\nix#1{}

\maketitle
\begin{abstract}
The purpose of this paper is to construct small-amplitude breather solutions for a nonlinear Klein-Gordon equation posed on a periodic metric graph 
via spatial dynamics and center manifold reduction.
The major difficulty occurs from the irregularity of the solutions. The persistence of the approximately constructed pulse solutions under higher order perturbations can be shown for two symmetric solutions.
\end{abstract}


\section{Introduction}

Differential operators on metric graphs with appropriate vertex conditions recently attracted a lot of interest, cf. \cite{KPost,Noja,Exner,Adami1}. 
For instance, they are used as simplified models 
for nano-technological objects with a similar geometric structure, such as nano-tubes or graphene. See the textbook \cite{Kuchment} for further motivations.

Korotyaev and Lobanov \cite{koro} found a unitarily equivalence between the periodic Schr\"o- dinger operator on a class of zigzag nanotubes and the direct sum of its corresponding Hamiltonians on a one-dimensional periodic metric graph 
consisting of rings and lines. 
Thus, they reduced the spectral problem on zigzag nanotubes to the spectral problem of periodic Schr\"odinger operators on one-dimensional graphs with necklace structure.
Recent works \cite{GPS16,PS17,MR3721872} studied nonlinear Schr\"odinger equations posed on this necklace graph.  
Particularly, the existence of small localized standing waves for frequencies lying below the linear spectrum of the associated stationary Schr\"odinger equation has been established in the work of Pelinovsky and Schneider \cite{PS17}.
By a variational approach, Pankov \cite{MR3721872} proved the existence of a non-small finite energy ground state solution.
The purpose of this article is to construct time-periodic, spatially localized solutions in nonlinear cubic Klein-Gordon equations.

From a mathematical point of view, the existence of these so called breather solutions is very rare. Breathers are an inherently nonlinear phenomenon.
In the spatially homogeneous situation, a wave equation known to admit small-amplitude breather solutions of pulse form is the Sine-Gordon equation.
However, these solutions do not persist under analytic perturbations. In particular, the Sine-Gordon equation is the only one of the form 
\begin{align*}
\partial_t^2u=\partial_x^2u-u-g(u)
\end{align*}
with $g(u)=\mathcal{O}(u^3)$ for $u\to0$ possessing breather solutions, \cite{Denzler}. 
In general, only the existence of generalized pulse solutions with small non-vanishing tails can be shown, cf. \cite{GS01, GS05, GS08}.
Our approach is motivated by the existence result of Blank, Chirilus-Bruckner, Lescarret, and Schneider \cite{BCLS11}. They considered a nonlinear Klein-Gordon equation 
\begin{align*}
s(x)\partial^2_tu-\partial_x^2u+q(x)u=u^3
\end{align*}
on the real line with specifically chosen, spatially periodic step functions $s$ and $q$. 
More recently, Hirsch and Reichel \cite{reichel} showed the existence of breather solutions of a semilinear wave equation with a periodically extended delta potential.
Crucial to their variational approach is that the spectrum of the corresponding wave operator is bounded away from zero.

The spectral picture necessary for the construction of breather solutions appears on the necklace graph in a natural way.
For a detailed spectral analysis we refer to \cite{molch,koro}.
The major difficulty will occur from the irregularity of the solutions caused by the imposed Kirchhoff boundary conditions, which lead to jumps of the first derivatives. 
As a consequence, the flow on the center manifold for the spatial dynamics formulation is no longer continuous as in \cite{BCLS11} with respect to the spatial evolution variable $x$. 

\subsection{Statement of the problem}

\begin{figure}
\begin{center}
\includegraphics[scale=1.3]{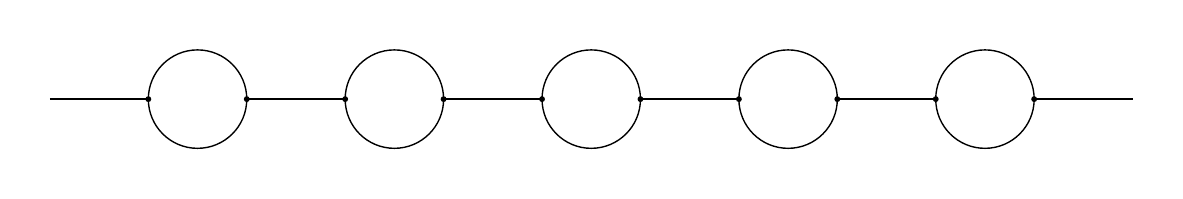}
\end{center}
\caption{The necklace graph is of the form $\Gamma=\oplus_{n\in\Z}\Gamma_n$ with $\Gamma_n=\Gamma^0_n\oplus\Gamma^+_n\oplus\Gamma^-_n$, 
where the $\Gamma^0_{n}$ are the horizontal links between the circles and the $\Gamma^{\pm}_{n}$ the upper and lower semicircles. 
The part $\Gamma^0_{n}$ is isometrically identified with the interval $I^0_{n}=[nP,nP+L]$ and the $\Gamma^{\pm}_{n}$ with the intervals $I^{\pm}_{n}=[nP+L,(n+1)P]$. 
The horizontal links are of lengths $L>0$, whereas the semicircles have length $\pi$. Hence, the periodicity of the graph is $P=L+\pi$.
For a function $u:\Gamma\to\C$, we denote the part on the interval $I_n^0$ with $u_n^0$ and the parts on the intervals $I_n^{\pm}$ with $u_n^{\pm}$.
}
\label{neckgr}
\end{figure}
We consider a cubic, nonlinear Klein-Gordon equation
\begin{align}\label{grundgl}
\partial_t^2u(t,x)=\partial_x^2u(t,x)-(\alpha+\varepsilon^2) u(t,x)+u(t,x)^3,\quad t\in\R,\, x\in\mathrm{int}\,\Gamma,
\end{align}
with a real-valued constant $\alpha$ and sufficiently small $\varepsilon>0$ on the periodic necklace graph $\Gamma$ from Figure \ref{neckgr}.
Throughout this paper we impose Kirchhoff boundary conditions at the vertex points $\{nP\}_{n\in\Z}$ and $\{nP+L\}_{n\in\Z}$, which consist of the continuity condition at the vertex points
\begin{align*}
u^0_n(nP+L)&=u^{\pm}_n(nP+L),\phantom{...}n\in\Z,\\
u^0_n((n+1)P)&=u^{\pm}_n((n+1)P),\,n\in\Z,
\end{align*}
and the conservation of the fluxes
\begin{align*}
\partial_xu_n^0(nP+L)&=\partial_xu_n^+(nP+L)+\partial_xu_n^-(nP+L),\phantom{.....}n\in\Z,\\
\partial_xu^0_{n+1}((n+1)P)&=\partial_xu_n^+((n+1)P)+\partial_xu_n^-((n+1)P),\,n\in\Z.
\end{align*}
It turns out that time-periodic, spatially localized solutions can be constructed within the invariant subspace of functions that are symmetric with respect to the semicircles.
In this case, the necklace graph $\Gamma$ can be identified with the real line equipped with a very singular periodic potential.

\subsection{Our main theorem}

Now we can state our main theorem.
\begin{thm}\label{main}
Let $L\in\{l\pi,\, l\in\N_{\mathrm{odd}}\}$ be the length of the horizontal links. 
For an odd integer $k$ and a sufficiently small $\varepsilon>0$ the nonlinear, cubic Klein-Gordon equation
\begin{align}\label{first}
\partial_t^2u(t,x)=\partial_x^2u(t,x)-\left(\frac{k^2}{4}+\varepsilon^2\right) u(t,x)+u(t,x)^3,\quad t\in\R,\, x\in\mathrm{int}\,\Gamma=\R\backslash(\pi\Z),
\end{align}
with Kirchhoff boundary conditions at the vertices possesses 
breather solutions of amplitude $\mathcal{O}(\varepsilon)$ and frequency $\omega=k/2$.
These solutions are symmetric in the upper and lower semicircles. 
Precisely, there exist functions $u:\R\times\R\to\R$ satisfying
\begin{itemize}
\item $u(t,x)=u(t+\frac{2\pi}{\omega},x)$ for all $t,x\in\R$,
\item $\mathrm{lim}_{\abs{x}\to\infty}u(t,x)e^{\beta\abs{x}}=0$ for all $t\in\R$ and a constant $\beta>0$. \fish
\end{itemize}
\end{thm}

\begin{rem}
The major challenge is the irregularity of the solutions due to Kirchhoff boundary conditions $(u^0)'(2n\pi)=2(u^+)'(2n\pi)$ and $2(u^0)'((2n+1)\pi)=(u^+)'((2n+1)\pi)$, $n\in\Z$, 
which leads to a non-autonomous system that makes it necessary to modify the persistence proof of the approximately constructed pulse under higher order perturbations. 
In contrast to the previous work the first derivative has jumps. As a consequence the flow on the center manifold is no longer continuous at the vertex points $x\in\pi\Z$. 
\end{rem}

\begin{rem}
Breathers are an inherently nonlinear phenomenon.
In principle, our method of proof allows to treat any odd nonlinearities, i.e.
\begin{align*}
N(u)=u^k,\quad k\in\N_{\mathrm{odd}}.
\end{align*}
\end{rem}

\subsection{Outline of the proof}

Using Fourier series expansion with respect to time $u(t,x)=\sum_{m}u_m(x)e^{im\omega t}$,
we transform the evolutionary problem (\ref{grundgl}) into countably many coupled second order ordinary differential equations for the Fourier coefficients
\begin{align}\label{spatial}
-m^2\omega^2u_m(x)=\partial_x^2 u_m(x)-(\alpha+\varepsilon^2) u_m(x)+(u*u*u)_m(x),\quad x\in\R,\, m\in\N_{\mathrm{odd}},
\end{align}
with new dynamic variable $x$, the so called spatial dynamics formulation.  
The cubic nonlinearity transforms into a discrete convolution.
Since we are interested in spatially localized solutions, i.e.
\begin{align*}
\mathrm{lim}_{\abs{x}\to\infty}u(t,x)=0,\quad t\in\R,
\end{align*}
we construct a homoclinic orbit to zero in the phase space of this infinite dimensional system (\ref{spatial}).
The key idea is to perform a center manifold reduction in order to reduce \eqref{spatial} to a finite dimensional system. 
However, because of the Kirchhoff boundary conditions, the system is non-autonomous and the first derivatives of the solutions have jumps and the flow on the center manifold is no longer continuous. 
Therefore, we apply a discrete version of the center manifold theorem to the family of time-$P$-maps.

We explain the core of our argumentation, which makes use of Floquet-Bloch theory. 
Linearizing \eqref{spatial} at the origin leads to the (decoupled) spectral problems
\begin{align}\label{spatial_lin}
-\partial_x^2 u_m(x)=(m^2\omega^2-\alpha) u_m(x)=\lambda_mu_m(x),\quad x\in\mathrm{int}\,\Gamma,\quad m\in\N_{\mathrm{odd}},
\end{align}
with Kirchhoff boundary conditions at the vertex points. 
Let $M(\lambda_m)$ denote the monodromy matrix of \eqref{spatial_lin}, which is the canonical fundamental matrix evaluated after one period of the system and conjugated 
to the linearizations of the time-$P$-maps at the origin. 
The complex number $\lambda_m$ corresponds to the spectrum of the negative Laplacian on the necklace graph if and only if the eigenvalues of $M(\lambda_m)$ (Floquet multipliers) lie on the complex unit circle.
Further, the number of Floquet multipliers on the complex unit circle determines the dimension of the center manifold.
Thus, we shall choose the constants $\omega$ and $\alpha$ such that $\lambda_1$ corresponds to the spectrum, 
whereas the positive numbers $\lambda_m$ for $3\leq m\in\N_{\mathrm{odd}}$ fall into spectral gaps of the negative Laplacian. 
Hence, the infinite dimensional spatial dynamics system \eqref{spatial} can then be reduced to a two-dimensional system on the center manifold.
Our method of construction heavily relies on the spectral properties of the linear system. In particular, the spectral gaps open linearly. 
In fact, these abstract center manifold constructions are related to a single ordinary differential equation
\begin{align*}
\partial_x^2u_c(x)=\varepsilon^2 u_c(x)-u_c^3(x),\quad x\in\mathrm{int}\,\Gamma,
\end{align*}
with imposed Kirchhoff boundary conditions on the graph, which will appear as the lowest order approximation of the dynamics on the center manifold.
The existence of a pulse solution has been established in \cite{PS17} via a detailed analysis of the stable and unstable manifold of the time-$P$-mapping.
In order to show persistence of this homoclinic orbit under higher order perturbations, we use reversibility and symmetry arguments.

The plan of the paper is as follows. 
In Section \ref{sec_spat} we introduce the spatial dynamics formulation and its symmetries.
The family of time-$P$-maps is investigated in Section \ref{subsec_timep}.
Section \ref{sec_spec} is dedicated to the linear spectral analysis on the periodic metric graph $\Gamma$. Moreover, we explain how to choose an adequate breather frequency 
and apply a discrete version of the center manifold theorem to the family of time-$P$-maps.
In Section \ref{sec_red} we relate these abstract center manifold constructions to a nonlinear cubic ordinary differential equation and find a homoclinic orbit, which persists under higher order perturbations.
Section \ref{discussion} contains a short discussion about arbitrary horizontal lengths of the necklace graph.
Finally, we have Appendices \ref{subsec_fb} and \ref{subsec_cm}, which contain a short introduction in Floquet-Bloch theory and discrete center manifold reductions.
\\

{\bf Funding.}  
This work was supported by the Deutsche Forschungsgemeinschaft DFG through the
Research Training Center GRK 1838 ``Spectral Theory and Dynamics of Quantum Systems". 

\section{Spatial dynamics formulation}\label{sec_spat}

The central purpose of this article is to find time-periodic, spatially localized solutions of the nonlinear Klein-Gordon equation
\begin{align}\label{first2}
\partial_t^2u(t,x)=\partial_x^2u(t,x)-(\alpha+\varepsilon^2) u(t,x)+u(t,x)^3,\quad t\in\R,\, x\in\mathrm{int}\,\Gamma,
\end{align}
with a real-valued constant $\alpha$ and $u(t)\in D(\partial_x^2|_{\Gamma})=\{C^0(\Gamma)\cap H^2(\mathrm{int}\,\Gamma)$: fulfilling Kirchhoff b.c. at $\partial\Gamma\}$. 
The Laplacian with domain $D(\partial_x^2|_{\Gamma})$ is self-adjoint, cf. \cite{Kuchment}.
Breather solutions can be constructed within the invariant subspace of symmetric functions with respect to the semi-circles. 
In this case, equation \eqref{first2} can be regarded as a Klein-Gordon equation on the real line equipped with a very singular periodic potential.
Searching for $\frac{2\pi}{\omega}$-periodic solutions
\begin{align*}
u(t,x)=u\left(t+\frac{2\pi}{\omega},x \right),\quad t,x\in\R,
\end{align*}  
Fourier series expansion leads to
\begin{align}
u(t,x)=\sum_{m\in\Z}u_m(x)e^{im\omega t}.
\end{align}
The real-valued constant $\omega$ has to be chosen suitably later on.
Thus, the evolutionary problem (\ref{first2}) transforms into countably many coupled second order ordinary differential equations 
\begin{align}\label{second}
-m^2\omega^2u_m(x)=\partial_x^2 u_m(x)-(\alpha+\varepsilon^2) u_m(x)+(u*u*u)_m(x),\quad m\in\Z,
\end{align}
where the cubic nonlinearity is given by a discrete convolution
\begin{align*}
((u*u*u)_m)_{m\in\Z}=\left(\sum_{n_1,n_2\in\Z}u_{m-n_1}u_{n_1-n_2}u_{n_2}\right)_{m\in\Z}.
\end{align*}
The dimension of the problem can be reduced by considering symmetries of the problem. Real-valued solutions satisfy $u_m=\overline{u_{-m}}$, $m\in\Z$. 
Moreover, the system is invariant under the transform $(t,u,u')\mapsto (-t,-u,-u')$, which leads to the condition $u_m=-u_{-m}$, $m\in\Z$. 
As an immediate consequence of the cubic nonlinearity, the space of solutions with $u_{2m}=0$, $m\in\Z$, is an invariant subspace. 
These conditions particularly lead to $\mathrm{Re}(u_m)=0$, $m\in\Z$. 
We prefer to replace $u_m$ by $iu_m$, where $iu_m$ satisfies the same equation with an opposite sign in front of the nonlinearity.
To conclude, we consider solutions in the invariant subspace
\begin{align*}
\hat{X}=\{(u_m)_{m\in\Z}\,|\,u_m\in\R, \, u_{2m}=0,\, u_m=u_{-m},\,m\in\Z\}
\end{align*}
of the system
\begin{align}\label{spatialsym}
-m^2\omega^2u_m(x)=\partial_x^2 u_m(x)-(\alpha+\varepsilon^2) u_m(x)-(u*u*u)_m(x),\quad m\in\Z_{\mathrm{odd}},
\end{align}
with Kirchhoff boundary conditions at the vertex points.

\section{Time-$P$-maps}\label{subsec_timep}

The first order system of \eqref{spatialsym} reads
\begin{align}\label{firstorder}
\partial_x\left(\begin{array}{c}u_m\\
u_m'\end{array}\right)=\left(\begin{array}{cc}0&1\\
-m^2\omega^2+\alpha&0\end{array}\right)\left(\begin{array}{c}u_m\\
u_m'\end{array}\right)+\left(\begin{array}{c}0\\
\varepsilon^2 u_m+(u*u*u)_m\end{array}\right),\, m\in\Z_{\mathrm{odd}},
\end{align}
with Kirchhoff boundary conditions at the vertex points. 
Interpreting the bifurcation parameter $\varepsilon$ as an independent variable, we treat the terms $\varepsilon^2 u_m$ in \eqref{firstorder} as nonlinear and use the denotation 
\begin{align}\label{form1}
\partial_xv_m=\Lambda_mv_m+N_m(\varepsilon,v),\, m\in\Z_{\mathrm{odd}}.
\end{align}
Denote by
\begin{align}\label{timeP}
v_{n,m}(\check{x})= \left(\begin{array}{c}u_m(\check{x}+nP;\check{x},\check{v})\\
u'^+_m(\check{x}+nP;\check{x},\check{v})\end{array}\right),\quad \check{x}\in[0,P),\,n\in\Z,
\end{align}
a solution $v_m$ at the vertex points $\check{x}+nP$ with initial conditions $v(\check{x})=\check{v}$ given at $\check{x}$.
We agree upon using right-hand sided derivatives at $\check{x}\in\{0,L\}$, since the Kirchhoff boundary conditions lead to jumps of the first derivative at the vertex points.
Now, the action of the time-$P$-mappings associated to \eqref{form1} can be written as
\begin{align}
T^Pv_{n,m}(\check{x}):=v_{n+1,m}(\check{x}),\quad n\in\Z,\,m\in\Z_{\mathrm{odd}}.
\end{align}
(The denotation time-$P$-map is chosen, since the spatial variable $\check{x}$ is the new dynamic variable.) 
The standard uniqueness theorem for second order ordinary differential equations with non-vanishing coefficient at the second derivative claims that if a solution and its derivative vanish at a point, 
it is identical zero. Therefore, the vector $v_{n,m}$ is well-defined on the invariant subspace of symmetric functions.
The linearizations at the origin of the associated discrete dynamical systems \eqref{timeP} are decoupled and 
given by monodromy matrices $M_{\check{x}}(m^2\omega^2-\alpha)$.
These matrices coincide with the canonical fundamental matrix of \eqref{form1} evaluated after one period of the problem and are conjugated to each other for $\check{x}\in[0,P)$, cf. Appendix \ref{subsec_fb}. 
For example, we explicitly compute for $\check{x}=0$,
\begin{align}\label{op_comm}
M_0(m^2\omega^2-\alpha)=\left(\begin{array}{cc} 1&0 \\
      0&2\end{array}\right)e^{\Lambda_m\pi}\left(\begin{array}{cc} 1&0 \\
      0&\frac{1}{2}\end{array}\right)e^{\Lambda_mL}.
\end{align}
In particular, we find
\begin{align}\label{form}
v_{n+1}=\Lambda_{\check{x}} v_n+\check{N}(\check{x},\varepsilon,v_n),\quad n\in\Z,
\end{align}
with a linear operator $\Lambda_{\check{x}}=\mathrm{diag}(M_{\check{x}}(m^2\omega^2-\alpha))_{m\in\Z_{\mathrm{odd}}}$ and a nonlinear map $\check{N}(\check{x},\varepsilon,v)$.
Since the matrices $M_{\check{x}}(m^2\omega^2-\alpha)$ are conjugated to each other, the eigenvalues of $\Lambda_{\check{x}}$ do not depend on $\check{x}$.
Our objective is to apply a discrete center manifold reduction to system \eqref{form}, cf. Appendix \ref{subsec_cm}. 

\section{Spectral situation and center manifold reduction}\label{sec_spec}

The discrete center manifold Theorem \ref{center} states that the number of eigenvalues of the linear operator $\Lambda$ lying on the unit circle is equal to the dimension of the center manifold.
The essential hypothesis \ref{spec} is the spectral separation of $\Lambda$, 
which requires a spectral gap around the unit circle.
Motivated by the important relation for the monodromy matrices, 
\begin{center}
$\abs{\mathrm{tr}(M(\lambda))}\leq 2$ $\Leftrightarrow$ two eigenvalues (Floquet multiplier) on the complex unit circle,
\end{center}
cf. Appendix \ref{subsec_fb}, we adjust the parameters $\omega$ and $\alpha$ in \eqref{spatialsym}, 
such that $\abs{\mathrm{tr}(M(\omega^2-\alpha))}\leq 2$ and $\abs{\mathrm{tr}(M(m^2\omega^2-\alpha))}> 2$ for any odd number $m\geq 3$. 
Hence, there will appear two Floquet multipliers on the unit circle, 
which lead to a family of two-dimensional center manifolds, cf. Subsections \ref{subsec_mon} and \ref{subsec_discred}. 
\begin{rem}
A detailed analysis of the spectrum of the Laplacian on the necklace graph is given in \cite{molch}. 
In contrast to periodic self-adjoint elliptic second order differential operator on the real line, there occurs a point spectrum. 
The eigenfunctions on the necklace graph are given by simple loop states, which are anti-symmetric with respect to the semicircles and vanish at the horizontal links, \cite{Kuchment, GPS16}. 
\end{rem}

\subsection{Trace of the monodromy matrix and choice of the breather frequency}\label{subsec_mon}

\begin{figure}
\begin{center}
\includegraphics[scale=0.8]{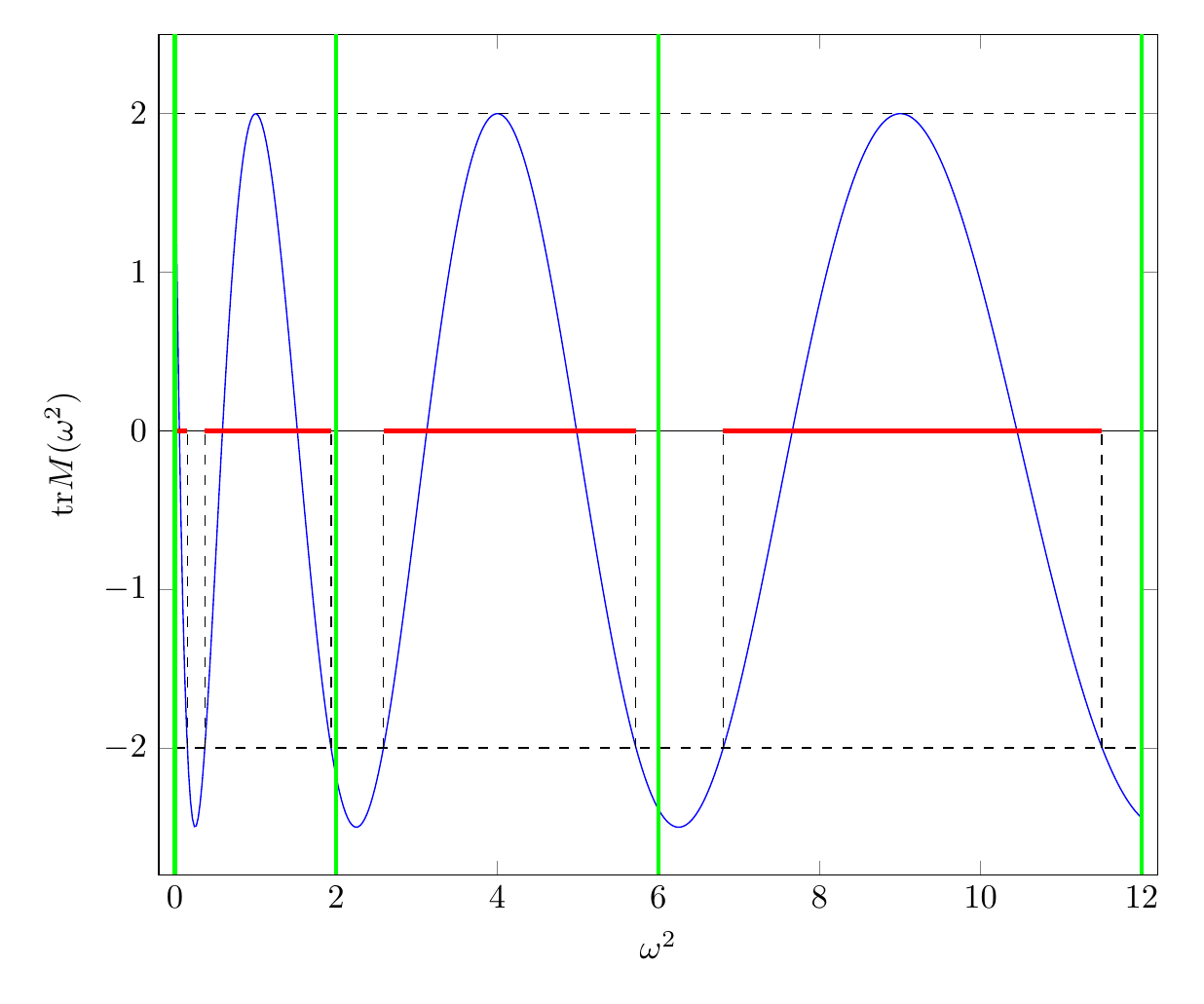}
\end{center}
\caption{
Trace of the monodromy matrix for $L=\pi$. The areas satisfying $\abs{\mathrm{tr}M(\omega^2)}>2$ lead to spectral gaps of the Laplacian on the graph, whereas the red lines correspond to its spectrum.
The points $\{(m^2-1)\omega^2:\,3\leq m\in\N_{\mathrm{odd}}\}$ (green lines) fall into spectral gaps and the point $\omega^2=0$ touches the spectrum.}
\label{figure_trace}
\end{figure}

\begin{figure}\label{hans}
\begin{center}
\includegraphics[scale=1]{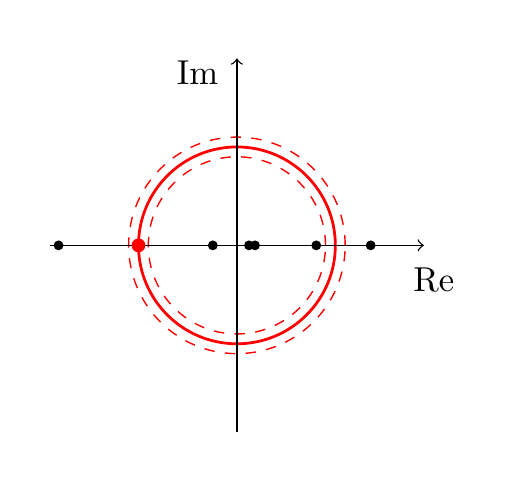}
\end{center}
\caption{We adjust the breather frequency such that there are two Floquet multiplier on the unit circle, which collide at $-1$, and the others are bounded uniformly away.}
\end{figure}

Let $L\in\{l\pi:\,l\in\N_{\mathrm{odd}}\}$. Taylor expansion of the exponential function in \eqref{op_comm} leads to
\begin{align}\label{11}
\mathrm{tr}M(\omega_m^2)=\frac{1}{4}(9\,\mathrm{cos}((L+\pi)\omega_m)-\mathrm{cos}((L-\pi)\omega_m)),
\end{align}
with $\omega_m^2=m^2\omega^2-\alpha$, cf. Figure \ref{figure_trace}.
The mapping $\omega\mapsto \mathrm{tr}M(\omega^2)$ is $1$-periodic. 
One of the major tasks is finding an adequate frequency $\omega$ and a constant $\alpha$, 
such that the operator $\Lambda=\mathrm{diag}(M(m^2\omega^2-\alpha))_{m\in\Z_{\mathrm{odd}}}$ has the property of spectral separation 
with precisely two Floquet multipliers lying on the unit circle, cf. Figure \ref{hans}. 
Choosing $\omega=k/2$, $k\in\N_{\mathrm{odd}}$ equal to an odd multiple of one half period, 
\begin{align}
\mathrm{tr}M(m^2\omega^2)=\begin{cases} -5/2, \quad\mathrm{for}\, \abs{m}\in\N_{\mathrm{odd}},\\
2,\phantom{-/5} \quad\mathrm{for}\, \abs{m}\in\N_{\mathrm{even}}.
\end{cases}
\end{align}
The fact that there are infinitely many Floquet multipliers on the unit circle prevents at a first view the application of the discrete center manifold theorem.
However, because of the symmetry of the spatial dynamics formulation, we only need to look for integers $m\in\Z_{\mathrm{odd}}$.
Varying the parameter $\alpha$, we can achieve $\abs{\mathrm{tr}M(\omega^2-\alpha)}=2$, whereas $\abs{\mathrm{tr}M(m^2\omega^2-\alpha)}>2$ for $\abs{m}\geq 3$.
In particular, we choose $\alpha=\omega^2$ in order to provoke the situation in the subsequent Theorem \ref{peli}.
Since the impact of $\alpha$ becomes smaller with increasing $m$, we find $\mathrm{tr}M(m^2\omega^2-\alpha)\to -5/2$ as $\N_{\mathrm{odd}}\ni \abs{m}\to\infty$.
As a consequence, we have two Floquet multipliers on the unit circle, which collide at $-1$.
Small perturbations of $\omega=k/2$, $k\in\N_{\mathrm{odd}}$, destroy this property.

\subsection{Discrete center manifold reduction}\label{subsec_discred}

The purpose of this subsection is to apply the discrete center manifold Theorem \ref{center} to the family of time-$P$-maps introduced in Section \ref{subsec_timep}.
We chose $\omega$ and $\alpha$ in Subsection \ref{subsec_mon} such that Floquet multipliers corresponding to the equations with index $m=\pm1$ collide at the unit circle, 
whereas the other ones are bounded uniformly away. 
Hence, the linear operator $\Lambda=(M(m^2\omega^2-\alpha))_{m\in\Z_{\mathrm{odd}}}$ has the property of spectral separation.
We comment on the other required properties of the linear and nonlinear part.
The linear operator $\Lambda$ maps the weighted sequence space $l^1_1(\Z,\R^2)$ into $l^1(\Z,\R^2)$ and Young's inequality sates
\begin{align*}
\norm{u*u*u}_{l^1_{\sigma}}\leq C\norm{u}_{l^1_{\sigma}}^3.
\end{align*}
The nonlinearity satisfies $\check{N}(\check{x})\in C^{\infty}(U\times V,l^2(\Z,\R^2))$ for a neighborhood $U\times V$ of $0$ in $\R\times l^1_1(\Z,\R^2)$ and
\begin{align*}
\check{N}(\check{x},0,0)=0, \quad \partial_3\check{N}(\check{x},0,0)=0.
\end{align*}
Let $\Pi_c$, respectively $\Pi_h$, denote the projections on the center and the hyperbolic subspace of the time-$P$-map $T^P$.  
We obtain
\begin{align}
\Pi_cv_n=(v_{n,1},v_{n,-1}),\quad \Pi_hv_n=(v_{n,3},v_{n,-3},v_{n,5},v_{n,-5},...).
\end{align}
Applying the discrete center manifold Theorem \ref{center}, it follows the existence of a reduction function $\Phi_{\check{x},\varepsilon}$ for $\varepsilon\in(-\varepsilon_0,\varepsilon_0)$, 
$\varepsilon_0>0$ sufficiently small, satisfying
\begin{align}
\Phi_{\check{x},\varepsilon}(\Pi_cv_n)=\Pi_hv_n.
\end{align}
Further, we remind of the symmetry restrictions in Section \ref{sec_spat},
\begin{align}
v_{n,-m}=v_{n,m}.
\end{align}
This leads to the family of two-dimensional reduced discrete systems
\begin{align}\label{equ_red}
v_{n+1,1}&=M(0)v_{n,1}+\tilde{N}_1(v_{n,1},\Phi_{\check{x},\varepsilon}(v_{n,1})),
\end{align}
with $\check{x}\in[0,P)$.

\section{Analysis of the reduced system}\label{sec_red}

\subsection{Relating the reduced discrete systems to an ordinary differential equation}

In order to analyze the reduced system \eqref{equ_red} on the center manifold we relate the abstract center manifold construction to an ordinary differential equation.
Let
\begin{align*}
v_{c,n}=\Pi_cv_n.
\end{align*}
The central projection $\Pi_c$ maps on the $u_1$-equation, respectively on the $(u_1,u'_1)$-part.
From the $P$-periodicity of the system we deduce
\begin{align}
v(\check{x}+(n+1)P;\check{x},\check{v})=v(\check{x}+P;\check{x},v(\check{x}+ nP;\check{x},\check{v})).
\end{align}
This leads to
\begin{align}
v_{c,n+1}=v_c(\check{x}+(n+1)P;\check{x},\check{v})=v_c(\check{x}+P;\check{x},v(\check{x}+nP;\check{x},\check{v})).
\end{align}
The projections $\Pi_{c}$ and $\Pi_h$ provide a decomposition of the Hilbert space $H$ into two invariant subspaces $H_c$ and $H_h$ 
and Theorem \ref{center} guarantees the existence of a reduction function $\Phi_{\check{x},\varepsilon}$. 
Since the nonlinearity does not possess quadratic terms, we deduce $v_h=\Phi_{\check{x},\varepsilon}(v_c)=\mathcal{O}(v_c^3)$. 
Thus, we derive
\begin{align}
v=(\Pi_cv)\oplus (\Pi_hv)=v_c\oplus \Phi_{\check{x},\varepsilon}(v_c)=v_c\oplus \mathcal{O}(v_c^3).
\end{align}
Inserting this relation leads to
\begin{align}
v_{c,n+1}&=v_c(\check{x}+P;\check{x},v(\check{x}+ nP;\check{x},\check{v}))\\ \nonumber
&=v_c(\check{x}+P;\check{x},v_c(\check{x}+nP;\check{x},\check{v})\oplus \Phi_{\check{x},\varepsilon}(v_c(\check{x}+nP;\check{x},\check{v})))\\ \nonumber
&=v_c(\check{x}+P;\check{x},v_c(\check{x}+nP;\check{x},\check{v})\oplus 0))+\mathcal{O}(v_c^5)\\ \nonumber
&=v_c(\check{x}+P;\check{x},v_{c,n}\oplus 0))+\mathcal{O}(v_c^5). 
\end{align}
Hence, in order to compute the flow on the center manifold up to $\mathcal{O}(v_c^5)$, it is sufficient to consider the discrete flow for $v_c$ on the center manifold with $v_h=0$. 
However, this discrete flow for any $\check{x}$ can be obtained by solving the ordinary differential equation for $u_1$ with $u_m=0$ 
for all $\abs{m}\geq 3$ and neglecting all terms of order $\mathcal{O}(u_1^4)$ and higher.
Hence, the cubic equation
\begin{align}
\partial_x^2u_1(x)&=\varepsilon^2 u_1(x)-u_1^3(x),
\end{align}
for $\varepsilon\in(0,\varepsilon_0)$ will appear as the lowest order approximation of the dynamics on the center manifold.

\subsection{Existence of the homoclinic orbit to zero on the center manifold}\label{subsec_ex}

We are exactly in the situation of \cite{PS17}, and recall the following
\begin{thm}\label{peli}
There are positive constants $\varepsilon_0$ and $C_0$, such that for every $\varepsilon\in(0,\varepsilon_0)$, the equation
\begin{align*}
\partial_x^2u=\varepsilon^2 u-u^3
\end{align*}
admits two non-trivial bound states $u\in D(\partial_x^2|_{\Gamma})$ (up to translational invariance) such that
\begin{align}
\norm{u}_{H^2(\Gamma)}\leq C_0\,\varepsilon.
\end{align}
One bound state satisfies 
\begin{align}\label{sym1}
u\left( x-\frac{L}{2}\right)=u\left( \frac{L}{2}-x\right) \,\mathrm{for}\,\mathrm{all}\, x\in\Gamma,
\end{align}
and the other one satisfies
\begin{align}\label{sym2}
u\left( x-\left( L+\frac{\pi}{2}\right)\right)=u\left( L+\frac{\pi}{2}-x\right) \,\mathrm{for}\,\mathrm{all}\, x\in\Gamma,
\end{align}
where $L$ is the length of the horizontal link and $\pi$ the length of the upper and lower semicircle. Moreover, the bound states obey the properties
\begin{itemize}
\item[i)] $u$ is symmetric in the upper and lower semicircles,
\item[ii)] $u(x)>0$ for every $x\in\Gamma$,
\item[iii)] $u(x)\to0$ as $\abs{x}\to\infty$ exponentially fast. \fish
\end{itemize}
\end{thm}
Hence, there exist two homoclinic orbits to the origin in the phase space.
It remains to prove persistence of these homoclinic solutions under higher order perturbations.  

\begin{figure}
     \setlength{\unitlength}{1cm}
 \begin{picture}(14,4)(0,0)
\put(0.4,-3){\includegraphics[width=7cm]{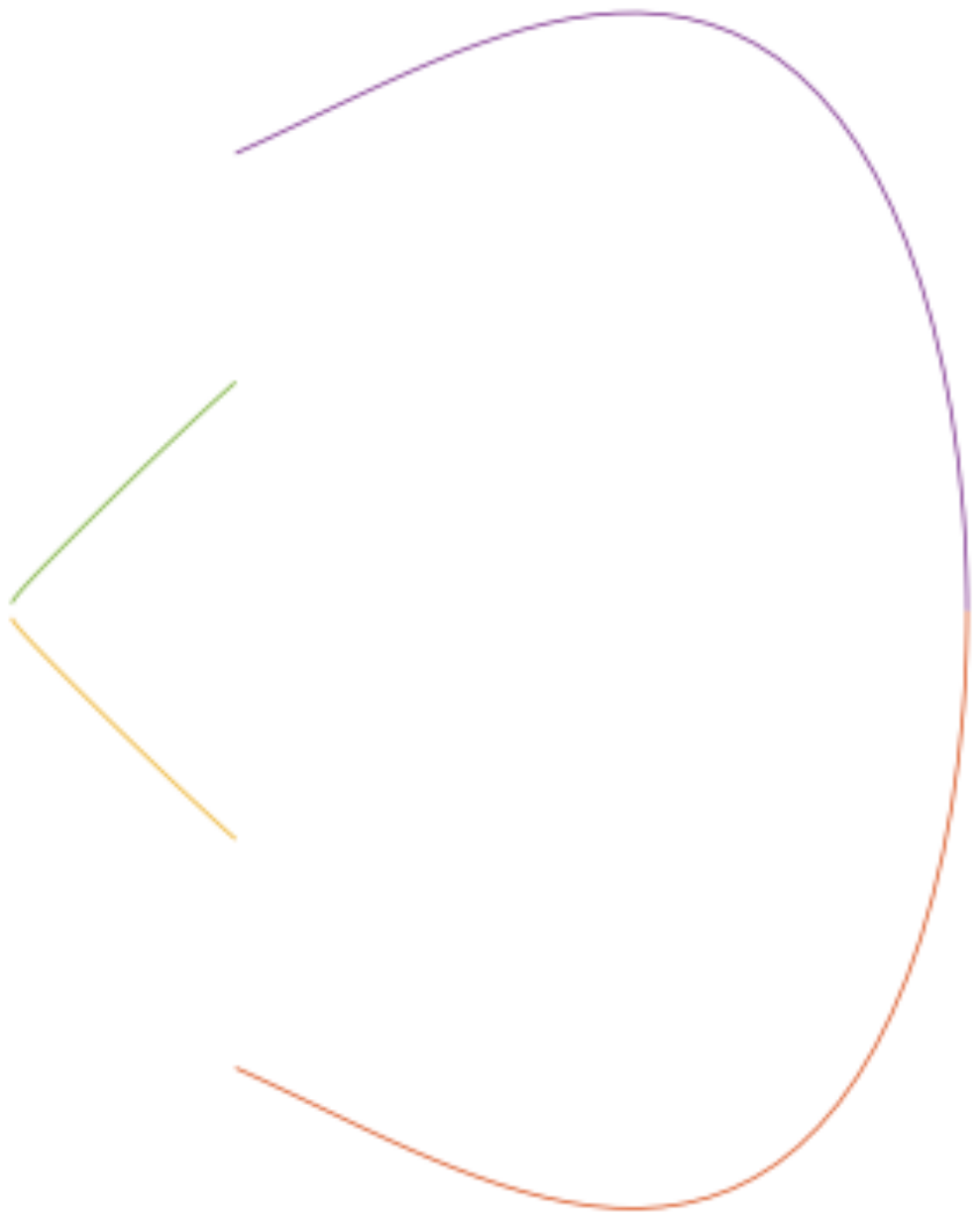}}
\put(1,2){\vector(1,0){6}}
\put(9.6,-0.8){\includegraphics[width=6cm]{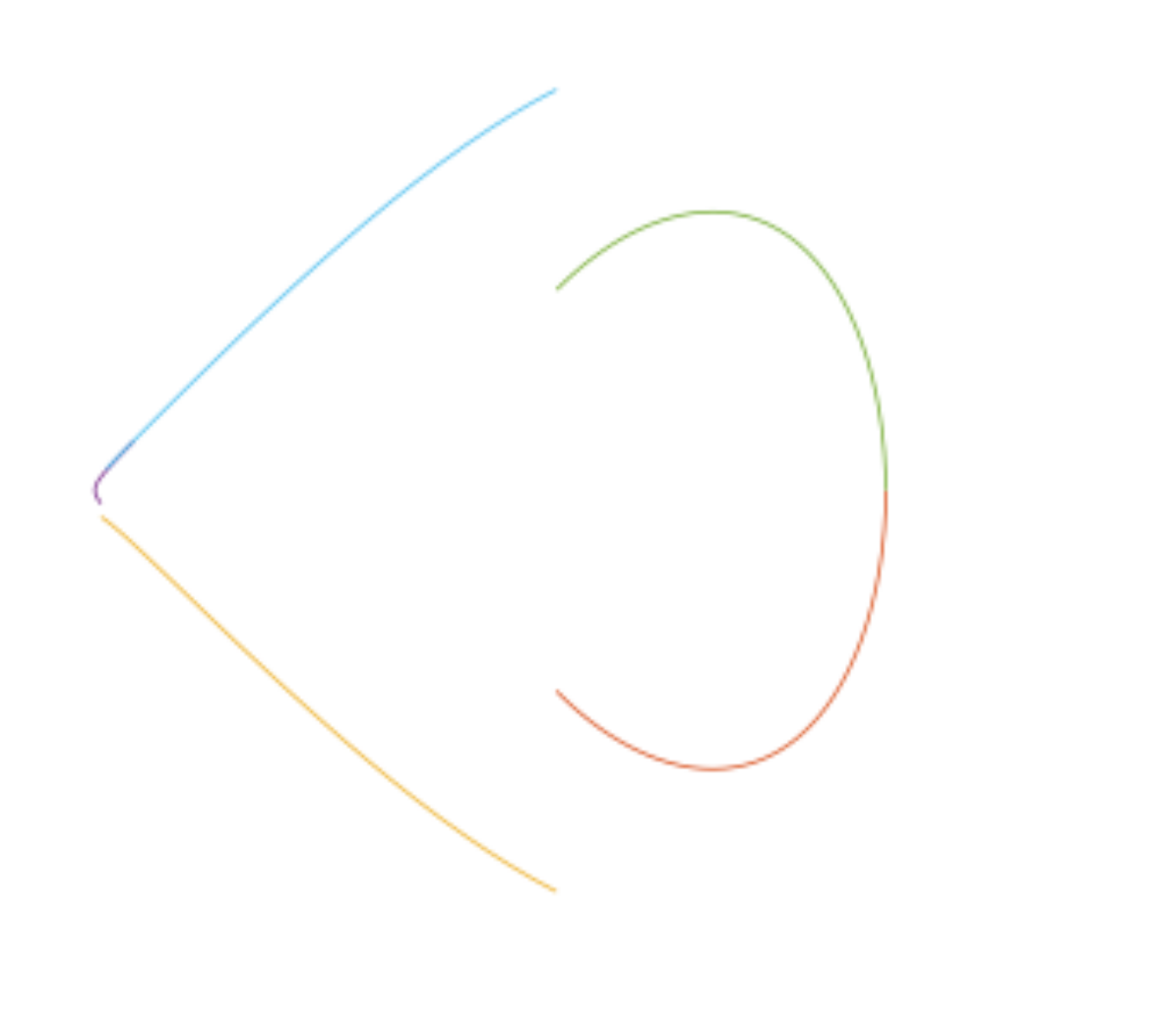}}
\put(6.5,1.5){$ u $}
\put(2,0){\vector(0,1){4}}
\put(9.6,3.5){$ u' $}
\put(9.1,2){\vector(1,0){6}}
\put(14.6,1.5){$ u $}
\put(10.1,0){\vector(0,1){4}}
\put(1.5,3.5){$ u' $}
\end{picture}
\vspace{0.7cm}
\caption{Sketch of homoclinics to the origin. The solution in the left panel corresponds to symmetry w.r.t. the midpoint of the horizontal link and in the right panel w.r.t. the midpoint of the semicircles.}
 \end{figure}

\subsection{Persistence of the homoclinic orbit under higher order perturbations}

The homoclinic orbit lies in the intersection of the stable and the unstable manifold. In general this intersection will break up, if higher order terms are added.
However, the situation is different in reversible systems. 
By proving a transversal intersection of the stable manifold with the fixed space of reversibility, we can construct the homoclinic solutions by reflecting the semi-orbit for $x\in(-\infty,x_0]$ at the $u$-axis.

First, the symmetries \eqref{sym1} and \eqref{sym2} are satisfied if and only if
\begin{align*}
\eqref{sym1}&\Leftrightarrow \partial_xu(L/2)=0,\\
\eqref{sym2}&\Leftrightarrow \partial_xu(L+\pi/2)=0,
\end{align*}
according to \cite{PS17}.
As a consequence, the homoclinic orbits 
intersect the $u$-axis transversally in the sense of smooth manifolds.
The smooth parts of length $\pi$, respectively $L$, are of order $\varepsilon^2$ in the $(u,u')$-plane. The perturbation will be of order $\varepsilon^3$ and cannot destroy the transversal intersection.
Second, the spatial dynamics system \eqref{spatialsym} 
is reversible, i.e. invariant under the mapping
\begin{align}\label{eqrev}
(x_0+x,u_m)\mapsto (x_0-x,u_m) 
\end{align}
for $x_0\in\{L/2,L+\pi/2\}$ due to the periodic structure of the graph and standard ordinary differential equation theory. 
Moreover, the corresponding time-$P$-maps admit a cut-off preserving reversibility, because they are derived from an even order explicit recurrence relation. We shall refer to \cite{james}, Section 5.2.
According to theorem \ref{rev}, the reduced system on the center manifold is also invariant under the mapping \eqref{eqrev}.
Therefore, as explained above, the approximative homoclinic solutions persist under higher order perturbations and exist in the full system, too.

\section{Discussion}\label{discussion}

Our previous argumentation heavily relies on the fact that the trace of the monodromy matrix is periodic. 
For general lengths $L$ this trace can be written as sum of $\mathrm{cos}$-terms, one of period $2\pi/(L+\pi)$ and the other one
of period $2\pi/(L-\pi)$, cf. equation \eqref{11}.
The sum of two periodic functions is periodic if and only if the ratio of the two periodicity constants is a rational number. 
Therefore, the trace is periodic if and only if
\begin{align*}
\frac{2\pi}{L+\pi}/\frac{2\pi}{L-\pi}=\frac{L-\pi}{L+\pi}\in\Q,
\end{align*}
which is the case for $L=l\pi$, $l\in\Q$. 
The existence of breather solutions can be established if the trace of the monodromy matrix evaluated at one half period has absolut value greater than two.
For instance, this is not fulfilled for an even integer $l$.
Moreover, we do not expect the existence of breathers for horizontal links of length $L=l\pi$, $l\notin\Q$.

\begin{rem}
Since we do not expect breather solutions if the ratio of the lengths is irrational, we predict that breather solutions will not persist under perturbations of the lengths.
\end{rem}

\appendix

\section{Floquet-Bloch theory}\label{subsec_fb}

To investigate the spectral problem
\begin{align}\label{allg}
-\partial_x^2u=\lambda u,\quad x\in\mathrm{int}\,\Gamma,
\end{align}
with constants $\lambda\in\R$ and Kirchhoff boundary conditions at the vertex points, 
we use tools from Floquet-Bloch theory, cf. \cite{RS, MR0404749, Kuchment}.
Let $u_1$ be the solution of \eqref{allg} with $u_1(0)=1$ and $u'^+_1(0)=0$ and let $u_2$ be the solution with $u_2(0)=0$ and $u'^+_2(0)=1$, where the index $+$ denotes the right-hand sided derivative. 
Consider the $2\times 2$ matrix
\begin{align}\label{monodef}
M_0(\lambda)=\left(\begin{array}{cc} u_1(P)&u_2(P) \\ u'^+_1(P) &u'^+_2(P)\end{array}\right),
\end{align}
which is a natural object, for if $v$ is a solution of (\ref{allg}), then
\begin{align}
\left(\begin{array}{c} v(P) \\ v'^+(P)\end{array}\right)=M_0(\lambda)\left(\begin{array}{c} v(0) \\ v'^+(0)\end{array}\right).
\end{align}
This means that the monodromy matrix $M_0(\lambda)$ is the fundamental matrix of the system of ordinary differential equations evaluated at the period of the system.

\begin{rem}\label{indep}
The monodromy matrices with varying evaluation points $\check{x}\in[0,P)$ are conjugated to each other.
Thus, their eigenvalues are independent of the evaluation points and so is $\mathrm{tr}(M_{\check{x}})=\mathrm{tr}(M)$ and $\mathrm{det}(M_{\check{x}})=\mathrm{det}(M)$.
\end{rem}

More insights gives the following

\begin{thm}[Floquet's Theorem]\label{Floquetthm}
There are linearly independent solutions $\psi_1,\psi_2$, such that either
\begin{itemize}
\item[i)] $\psi_1(x)=e^{m_1x}p_1(x)$ and $\psi_2(x)=e^{m_2x}p_2(x)$,
or
\item[ii)] $\psi_1(x)=e^{mx}p_1(x)$ and $\psi_2(x)=e^{mx}(xp_1(x)+p_2(x))$,
\end{itemize}
with constants $m_1,m_2,m\in\C$ and $P$-periodic functions $p_1,p_2$. \fish
\end{thm}

In other words, Floquet's theorem shows that the fundamental matrix $\Phi(x)$ with $\Phi(0)=I$ can be written as
\begin{align*}
\Phi(x)=Q(x)e^{xN},
\end{align*}
with $Q(x+P)=Q(x)$ and a matrix $N$ independent of $x$, which is similar to a diagonal matrix in the case i) and has a Jordanblock in the case ii).
We want to emphasize the simple connection $M=e^{PN}.$ between the monodromy matrix $M$ defined in \eqref{monodef} and the matrix $N$.
Therefore, we deduce
\begin{align*}
\mu_i=e^{m_iP},
\end{align*}
where $m_i$ are the constants of Theorem \ref{Floquetthm} and $\mu_i$ denotes the eigenvalues of the monodromy matrix. 
The monodromy matrix is known to have determinant $1$, which implies that its eigenvalues are $\mu$ and $\mu^{-1}$ and $\mathrm{tr}(M)=\mu+\mu^{-1}$. 
We can distinguish the following cases:
\begin{itemize}
\item[1)] $\mathrm{tr}(M)>2:$ The eigenvalues $\mu_1\neq\mu_2$ are positive, real numbers not equal to $1$ and the linearly independent solutions are exponentially growing/decaying and of the form
\begin{align*}
\psi_{1,2}(x)=e^{\pm mx}p_{1,2}(x)
\end{align*}
with a positive constant $m$.
\item[2)] $\mathrm{tr}(M)<-2:$ The eigenvalues $\mu_1\neq\mu_2$ are negative, real numbers not equal to $-1$ and the linearly independent solutions are exponentially growing/decaying and of the form
\begin{align*}
\psi_{1,2}(x)=e^{(\pm m+i\pi/P)x}p_{1,2}(x)
\end{align*}
with a positive constant $m$.
\item[3)] $-2<\mathrm{tr}(M)<2:$ The eigenvalues $\mu_1\neq\mu_2$ lie on the complex unit circle away from $\{\pm 1\}$. The eigenfunctions are uniformly bounded and 
\begin{align*}
\psi_{1,2}(x)=e^{\pm ilx}p_{1,2}(x)
\end{align*}
with a real constant $l$.
\item[4)] $\mathrm{tr}(M)=\pm 2:$ In this case the eigenvalues are equal to $\{\pm 1\}$. The second part of theorem \ref{Floquetthm} applies if and only if $M$ is similar to the Jordanblock
\begin{align*}
\left(\begin{array}{cc} \pm 1&1 \\
      0&\pm 1\end{array}\right),
\end{align*}
and this is the case, if and only if $\mathrm{tr}(M)$ has a turning point at $\pm 2$. Otherwise, part i) applies and we have two periodic eigenfunctions in the case $\mu_1=\mu_2=1$, respectively semi-periodic for $\mu_1=\mu_2=-1$.
\end{itemize}

\noindent To sum up, we have the following equivalences
\begin{center}
$\abs{\mathrm{tr}(M(\lambda))}\leq 2$ $\Leftrightarrow$ Floquet multiplier on the complex unit circle $\Leftrightarrow$ $\lambda\in\sigma(-\partial_x^2|_{\Gamma})$\\
$\abs{\mathrm{tr}(M(\lambda))}> 2$ $\Leftrightarrow$  Floquet multiplier off the complex unit circle $\Leftrightarrow$ $\lambda\notin\sigma(-\partial_x^2|_{\Gamma})$
\end{center}

\section{The discrete center manifold theorem}\label{subsec_cm}

For the reader's convenience we recall a discrete version of the center manifold theorem and refer to \cite{james}.
First, we describe the general framework, in which the center manifold reduction applies. Let $H$ be a Hilbert space and consider a closed linear operator $\Lambda:D\subset H\to H$.
We equip $D$ with the scalar product $\langle u,v\rangle_D=\langle \Lambda u,\Lambda v\rangle_H+\langle u,v\rangle_H$, which leads to the Hilbert space $D$ continuously embedded in $H$.
Further, denote by $U\times V$ a neighborhood of $0$ in $\R\times D$ and assume that the nonlinear map $N\in C^k(U\times V,H)$ for at least $k\geq 2$ satisfies
\begin{align*}
N(0,0)=0, \quad D_uN(0,0)=0.
\end{align*}
We look for sequences $(y_n)_{n\in\Z}$ in $V$ satisfying
\begin{align}\label{dis}
y_{n+1}=\Lambda y_n+N(\varepsilon,y_n), \quad\mathrm{in}\,H,\quad \forall n\in\Z,
\end{align}
with a constant $\varepsilon$ independent of $n$.

\begin{rem}
The condition $N(0,0)=0$ means that $0$ is an equilibrium of the discrete equation, and the condition $D_uN(0,0)=0$ then shows that $\Lambda$ is the linearization of the vector field about $0$, so that $N$ represents the nonlinear terms, which are of the order $\mathcal{O}(\norm{y}_H^2)$. 
\end{rem}

\begin{hyp}\label{spec}
The operator $\Lambda$ has the property of spectral separation, which means that its spectrum $\sigma(\Lambda)$ splits in the following way
\begin{align*}
\sigma(\Lambda)=\sigma_s\cup\sigma_c\cup\sigma_u,
\end{align*} 
where $\sigma_s=\{z\in\C\,:\,\abs{z}<1\}$, $\sigma_c=\{z\in\C\,:\,\abs{z}=1\}$ and $\sigma_u=\{z\in\C\,:\,\abs{z}>1\}$. We further assume 
$\mathrm{sup}_{z\in\sigma_s}\abs{z}<1$ and $\mathrm{inf}_{z\in\sigma_u}\abs{z}>1$.
\fish
\end{hyp}

Hence, the hyperbolic part $\sigma_s\cup\sigma_u$ of the system has nonzero distance to the center part, i.e. there is a spectral gap around the unit circle, which allows us to define spectral projections:
\begin{align*}
\Pi_c&=\frac{1}{2\pi i}\int_{C(R)}(\mu-\Lambda)^{-1}d\mu-\frac{1}{2\pi i}\int_{C(r)}(\mu-\Lambda)^{-1}d\mu,\\
\Pi_h&=Id-\Pi_c,
\end{align*}
where $C(r)$ denotes the circle with center in zero and radius $r$ and 
\begin{align*}
\mathrm{sup}_{z\in\sigma_s}\abs{z}<r<1<R<\mathrm{inf}_{z\in\sigma_u}\abs{z}.
\end{align*}
We introduce some notation for the center space $H_c=\Pi_c H$, as well as the hyperbolic projection $\Pi_h=I_H-\Pi_c$ and $H_h=\Pi_hH$, $D_h=\Pi_hD$. 
The projections provide a decomposition of $H$ into the two invariant subspaces $H_c$ and $H_h$.

\begin{thm}[Discrete center manifold theorem]\label{center}
Under Hypothesis \ref{spec} there exists a neighborhood $I\times\Omega$ of $0$ in $\R\times D$ and a map $\Phi_{\varepsilon}\in C^k_b(I\times H_c,D_h)$ such that for all $\varepsilon\in I$ the manifold
\begin{align}
M_{\varepsilon}=\{y\in D: \, y=x+\Phi_{\varepsilon}(x), \,x\in H_c\}
\end{align}
has the following properties
\begin{itemize}
\item[i)] $M_{\varepsilon}$ is locally invariant, i.e. if $y\in M_{\varepsilon}\cap\Omega$, then $\Lambda y+N_{\varepsilon}(y)\in M_{\varepsilon}$.
\item[ii)] If $(y_n)_{n\in \Z}\subset\Omega$ is a solution of (\ref{dis}), then $y_n\in M_{\varepsilon}$ for all $n$ and the recurrence relation
\begin{align}\label{res}
y_{n+1}^c=f_{\varepsilon}(y_n^c),\quad\forall n\in\Z,
\end{align}
is satisfied in $H_c$, where the function $f\in C^k(I\times(H_c\cap\Omega),H_c)$ is defined by
\begin{align}
f_{\varepsilon}=\Pi_c(\Lambda+N_{\varepsilon})\circ(I+\Phi_{\varepsilon}).
\end{align}
\item[iii)] Conversely, if $(y_n^c)_{n\in\Z}\subset\Omega$ is a solution of (\ref{res}), then
\begin{align}
y_n=y_n^c+\Phi_{\varepsilon}(y_n^c), \quad n\in\Z,
\end{align}
satisfies (\ref{dis}). \fish
\end{itemize}
\end{thm}

The manifold $M$ is called a local center manifold and the map $\Phi$ is referred to as reduction function.
This theorem allows us to reduce the local study of the discrete equation (\ref{dis}) to that of the recurrence relation (\ref{res}) on the subspace $H_c$, which is particularly interesting when $H_c$ is finite dimensional. 

We finish this section with a reduction result preserving reversibility. 
Let \eqref{dis} be reversible with respect to a symmetry $R\in\mathcal{L}(D)$, i.e. if $u_n$ is a solution, then $Ru_{-n}$ is also a solution.

\begin{thm}\label{rev}
Assume additionally to the assumptions in Theorem \ref{spec} that $\Lambda$ admits a cut-off preserving a reversibility symmetry $R$, see \cite{james} Definition 2.
Then, the reduced mapping is reversible and one has
\begin{align*}
R\circ\Pi_c=\Pi_c\circ R,\quad R\circ\Phi_{\varepsilon}=\Phi_{\varepsilon}\circ R. \rtimes
\end{align*}
\end{thm}

\textbf{Acknowledgement.}
I would like to thank my supervisor Guido Schneider for all his support and encouragement during the realization of this paper.

\end{document}